\documentclass{amsproc}
\usepackage{amssymb}
\usepackage{graphicx}
\usepackage{bm}
\usepackage[nameinlink]{cleveref}
\usepackage{MnSymbol}
\usepackage{multirow}
\usepackage{makecell}
\usepackage{float}
\usepackage{color}

\newcommand{\bfPW}{\bm P\bm W}
\newcommand{\bfg}{\bm g}

\newcommand{\bfA}{\bm A}
\newcommand{\bfB}{\bm B}
\newcommand{\bfa}{\bm a}
\newcommand{\bfb}{\bm b}
\newcommand{\bfd}{\bm d}
\newcommand{\bfC}{\bm C}
\newcommand{\bfE}{\bm E}
\newcommand{\bfF}{\bm F}
\newcommand{\bfL}{\bm L}
\newcommand{\bfH}{\bm H}
\newcommand{\bfV}{\bm V}
\newcommand{\bfu}{\bm u}
\newcommand{\bfq}{\bm q}
\newcommand{\bfr}{\bm r}
\newcommand{\bft}{\bm t}
\newcommand{\bfv}{\bm v}
\newcommand{\bfnu}{\bm \nu}
\newcommand{\bfxi}{\bm \xi}
\newcommand{\bftau}{\bm \tau}
\newcommand{\bfchi}{\bm \chi}
\newcommand{\bfx}{\bm x}
\newcommand{\bfz}{\bm z}
\newcommand{\bfp}{\bm p}
\newcommand{\bfw}{\bm w}
\newcommand{\bfW}{\bm W}
\newcommand{\cT}{\mathcal{T}}
\newcommand{\cF}{\mathcal{F}}
\DeclareMathOperator{\curl}{\nabla\times}
\DeclareMathOperator{\divv}{\nabla\cdot}
\newcommand{\jmp}[1]{\lsem #1{} \rsem }   
\newcommand{\avg}[1]{\left\{\!\!\left\{#1\right\}\!\!\right\}}

\newtheorem{theorem}{Theorem}[section]

\theoremstyle{definition}

\theoremstyle{remark}

\numberwithin{equation}{section}
\newcommand{\rev}[1]{#1}
\begin{document}

\title{Finite Element Methods for Maxwell's Equations}


\author{Peter Monk}
\address{Department of Mathematical Sciences, University of Delaware, Newark DE 19716, USA}
\email{monk@udel.edu}
\thanks{}

\author{Yangwen Zhang}
\address{Department of Mathematical Sciences, University of Delaware, Newark DE 19716, USA}
\curraddr{}
\email{ywzhangf@udel.edu}
\thanks{}

\subjclass[2000]{Primary }

\date{}

\begin{abstract}
We survey finite element methods for approximating the time harmonic Maxwell equations.  We concentrate on
comparing error estimates for problems
with spatially varying coefficients.  For the conforming edge finite element methods, such estimates allow, at least, piecewise smooth coefficients.  But for Discontinuous Galerkin (DG) methods,
the state of the art of error analysis is less advanced (we consider three DG families of methods: Interior Penalty type, Hybridizable DG, and Trefftz type methods).  Nevertheless, DG methods offer
significant potential advantages compared to conforming methods.
\end{abstract}

\maketitle

\section{Introduction}
Maxwell's equations govern the propagation of electromagnetic waves, and so applications are
widespread in science and technology.  When the
material through which the waves propagate is
inhomogeneous (i.e. has spatially varying electromagnetic parameters) or is anisotropic,
it is common to discretize the Maxwell system directly (rather than via integral equations).
In this paper we shall survey finite element 
methods for approximating the solution of
frequency domain (or time-harmonic) electromagnetic wave propagation.  The main emphasis of the paper is on discontinuous 
Galerkin (DG) methods, but we shall first recall
a brief history of conforming methods and mention the important input that mathematics has 
provided in their development.

An obvious question is: Why use DG methods? The general answer is flexibility.  It is easy to change the local discretization spaces (for example, polynomial degree) from element to element, and typically hanging nodes are permissible.  More recently, the ability to use general polyhedral mesh elements has emerged as an attractive 
feature~\cite{houston17}. 
Other attractive features of DG methods include their element-wise conservation properties. Especially, comparing with $\bm H({\rm curl};\Omega)$-conforming edge element methods, it is much easier to implement a 
DG method when using hp-adaptivity for Maxwell's equations.

The usual objection
to DG methods is that, for a given mesh and polynomial degree, they have more degrees of freedom than conforming methods.  But this problem can be ameliorated by using a hybridizable DG
method, or a Trefftz method.

We shall only consider the simplest boundary value problem for Maxwell's equations.  To describe this problem, let $\Omega\subset \mathbb{R}^3$ denote a bounded connected Lipschitz polyhedral domain with connected boundary $\Gamma:=\partial\Omega$, and let $\bfnu$ denote the
unit outward normal.  We suppose that the material in $\Omega$ has relative magnetic permeability $\mu_r=1$ (i.e. magnetic effects are absent), 
but that the 
relative permittivity $\epsilon_r$ is a function
of position and for simplicity is piecewise
in $W^{1,\infty}$ (see \cite{Ass18}).  In general
$\Re(\epsilon_r)$ may be positive (for example
in a dielectric) or negative (as for certain metals like gold at optical frequencies).  However the conductivity $\Im(\epsilon_r)\geq 0$ and $\Im(\epsilon_r)>0$ when $\Re(\epsilon_r)<0$.  For simplicity we shall assume $\Re(\epsilon_r)$ is strictly positive in accordance with most of the papers we shall reference (where it is assumed that the wavelength is such that metals can be modelled as impenetrable regions). When $\Re(\epsilon_r)<0$ and $\Im(\epsilon_r)>0$, error estimates are usually a simple matter. This is not true if $\Re(\epsilon_r)<0$ and $\Im(\epsilon_r)=0$, and for sign changing coefficients with zero conductivity, see for example
~\cite{ASBB14}.

In our upcoming discussion vector quantities
or spaces of vector functions will be in bold
typeface. 

Suppose the electromagnetic field has angular frequency $\omega$. We define the wave number
$\kappa=\omega/c$ where $c$ is the speed of light
in vacuum. Then the electric field $\bfE$ satisfies the following Maxwell
system written in second order form:
\begin{subequations}\label{Maxwell}
\begin{align}
\curl\curl\bfE-\kappa^2\epsilon_r
\bfE&=\bfF\mbox{ in }\Omega,\label{Maxwell-curl}\\
\bfnu\times\bfE&=0\mbox{ on }\Gamma:=\partial\Omega.\label{Maxwell-PEC}
\end{align}
\end{subequations}
Here $\bfF$ describes imposed currents, and for simplicity we shall assume $\bfF\in \bfL^2(\Omega)$ and $\divv\bfF=0$.  The boundary
condition \cref{Maxwell-PEC}  is termed the
perfect electrical conductor or PEC boundary
condition and is appropriate for metal at microwave frequencies.

In order to write down variational formulations of \cref{Maxwell}, we need some notation.
Let $\Lambda\subset\mathbb{R}^3$ denote a bounded Lipschitz domain.  Then the standard function space for analysing \cref{Maxwell} is  
\[
\bfH({\rm curl};\Lambda):=\{\bfu\in \bfL^2(\Lambda)\;|\;\curl\bfu\in \bfL^2(\Lambda)\}
\]
in the case when $\Lambda=\Omega$. Taking into account the PEC boundary condition \cref{Maxwell-PEC} the appropriate subspace for the solution is
\[
\bfH_0({\rm curl};\Lambda):=\{\bfu\in \bfH({\rm curl};\Lambda)\;|\;\bfnu\times\bfu=0\mbox{ on }\Gamma\},
\]
with norm
\[
\Vert \bfu\Vert_{{\rm curl},\Lambda}^2=\Vert\curl \bfu\Vert_{\Lambda}^2
+\Vert \bfu\Vert_{\Lambda}^2
\]
where $\Vert\cdot\Vert_\Lambda$ denotes the $\bfL^2$ norm on $\Lambda$.  

The main tool to derive variational formulations is the following integration by parts identity
\begin{equation}
\int_\Lambda\curl \bfA\cdot\bfB\,dV=\int_\Lambda\bfA\cdot\curl\bfB\,dV+
\int_{\partial \Lambda}\bfnu\times\bfA\cdot\bfB_T\,dA,\label{intp}
\end{equation}
where $\bfB_T=\bfnu\times(\bfB\times\bfnu)$.  This identity is valid for functions $\bfA,\bfB\in \bfH({\rm curl};\Lambda)$ and for Lipschitz domains~\cite{Buffa-Ciarlet01,KHBook}. We also use the notation
\[
(\bfu,\bfv)_\Lambda:=\int_\Lambda\bfu\cdot\overline{\bfv}\,dV,
\]
where the over-bar denotes complex conjugation.

Proceeding formally we can derive a weak form for \cref{Maxwell} by multiplying the first equation by the complex conjugate of a
smooth test vector $\bfxi\in \bfC^{\infty}_0(\Omega)$ and using the integration by parts identity \cref{intp} with $\Lambda=\Omega$, to obtain
\begin{equation}
\int_\Omega\curl \curl\bfE\cdot\overline{\bfxi}\,dV=\int_\Omega\curl \bfE\cdot\curl\overline{\bfxi}\,dV+
\int_{\Gamma}\bfnu\times\bfE\cdot\overline{\bfxi}_T\,dA,
\label{maxbyparts}
\end{equation}
where the boundary term vanishes because $\bfxi_T$ vanishes on $\partial\Omega$.
Using this equality, we derive the variational problem of finding $\bfE\in \bfH_0({\rm curl};\Omega)$ such that
\begin{equation}
(\curl\bfE,\curl\bfxi)_\Omega-\kappa^2(\epsilon_r\bfE,\bfxi)_\Omega=(\bfF,\bfxi)_\Omega
\mbox{ for all }\bfxi\in \bfH_0({\rm curl};\Omega).\label{maxweak}
\end{equation}

The above problem (\ref{maxweak}) reveals two issues that complicate the analysis of Maxwell's equations:
\begin{enumerate}
    \item Only the curl of $\bfE$ appears in the derivative term in this equation. There is no explicit control over the divergence. But we note that by choosing
    \[
    \xi\in H^1_0(\Omega):=\{u\in L^2(\Omega)\;|\; \nabla u\in \bfL^2(\Omega) \mbox{ and }u=0\mbox{ on }\Gamma\},
    \]
    we have $\nabla\xi\in \bfH_0({\rm curl};\Omega) $. Using $\bfxi=\nabla\xi$ in
    \cref{maxweak} we obtain, using our assumption that $\divv\bfF=0$,
    \begin{equation}
    (\epsilon_r\bfE,\nabla\xi)_\Omega=0 \mbox{ for all }\xi\in H^1_0(\Omega).
    \label{divfree}
    \end{equation}
    Thus $\divv(\epsilon_r\bfE)=0$.   This suggest that the correct space for the analysis of \cref{maxweak} is then
    \[
    \bfV:=\bfH_0({\rm curl};\Omega)\cap \bfH({\rm div}_{\epsilon_r}^0,\Omega),
    \]
    where
    \[
    \bfH({\rm div}_{\epsilon_r}^0;\Omega):=\{\bfu\in \bfL^2(\Omega)\;|\;\divv(\epsilon_r\bfu)=0\mbox{ in }\Omega\}.
    \]
    Note that this implies that $\bfE$ generally has a discontinuous normal component across surfaces where $\epsilon_r$ is discontinuous.
    \item Even using $\rev{\bfV}$ in place of $\bfH_0({\rm curl};\Omega)$ in \cref{maxweak}, the sesquilinear form on the left hand side of \cref{maxweak} is not coercive.  However, it does satisfy a G\aa{}rding inequality, and then the compact embedding of $\rev{\bfV}$ into $\bfL^2(\Omega)$~\cite{Leis88} shows that \cref{maxweak} has a unique solution
    provided $\kappa$ is not an eigenvalue for the curl-curl operator~\cite{Monk03}.  We assume this restriction on $\kappa$ from now on.
\end{enumerate}
For more information about the mathematical theory of Maxwell's equations see~\cite{Monk03,KHBook,Ass18}.

The outline of the paper is as follows. In the next section we shall briefly discuss conforming discretizations of Maxwell's
equations before turning to the main subject of this paper: DG methods.  There are many versions of DG methods, and we have chosen just three examples to discuss here. In Section~\ref{IPDG} we discuss one of the first DG methods for Maxwell's equations, the Interior Penalty Discontinuous Galerkin Method (IPDG).  We then move, in Section \ref{HDG} to an alternative approach that allows for new solution strategies, the Hybridizable Discontinuous Galerkin (HDG) method, before continuing with a Trefftz type method, the Ultra Weak Variational Formulation (UWVF) in Section~\ref{UWVF}.  We end with some conclusions.

\section{Finite Elements and Maxwell's Equations}\label{FEMax}
In this section we shall discuss the standard conforming approximation to
Maxwell's equations.  The most direct approach would be to construct
a conforming approximation in $\rev{\bfV}$, but this brings certain dangers that 
we shall 
mention in Section \ref{stabh1}. Instead we prefer to work in $\bfH_0({\rm curl};\Omega)$ (see Section~\ref{edgel}).  

Although discontinuous Galerkin methods can be used on rather general meshes, we restrict ourselves in this paper to the simplest case, a tetrahedral mesh.  Let $\{\cT_h\}_{h>0}$ denote a family of meshes of shape regular tetrahedra that cover $\Omega$ and satisfy the usual finite element meshing constraints~\cite{Ciarlet78}. The parameter $h$ denotes the maximum diameter
of the elements in $\cT_h$.

\subsection{Stabilized $H^1$ elements}\label{stabh1}
Let us suppose in this subsection that $\epsilon_r$ is constant.  We could then
define a finite element subspace of $\rev{\bfV}$ of polynomials of degree $k$ by
\[
\bfV_h^{(1)}:=\{\bfu\in \rev{\bfV}\;|\; \bfu|_K\in \mathbb{P}_k^3\mbox{ for all }K\in \cT_h\},
\]
where $\mathbb{P}_k$ denotes the space of polynomials of total degree at most $k$ in three variables.
\rev{If $\bfu_h\in\bfV_h^{(1)}$, the} requirement that $\bfu_h\in \rev{\bfV}$ \rev{implies it must have a well} defined divergence and curl. This
necessarily requires that $\bfu_h$ is continuous and hence $\bfV_h^{(1)}\subset \bfH^1(\Omega)$ (both normal and tangential components of the field must be continuous across faces in the mesh).  Thus we can construct $\bfV_h^{(1)}$ from three copies of the standard scalar continuous piecewise $k$-degree  finite element spaces.  This is a major
practical attraction.

Imposition of the divergence free constraint suggests to modify \cref{maxweak} to include a penalty term for the divergence~\cite{Hazard-Lenoir96}, so choosing $\gamma>0$ we seek $\bfE_h\in \bfV_h^{(1)}$ such that
\begin{equation}
(\curl\bfE_h,\curl\bfxi)_\Omega+\gamma(\divv\bfE_h,\divv\bfxi)_{\Omega}-\kappa^2(\epsilon_r\bfE_h,\bfxi)_\Omega=(\bfF,\bfxi)_\Omega
    \label{maxweak1}
\end{equation}
for all $\bfxi\in \bfV_h^{(1)}$.  

Unfortunately this formulation has subtle problems
that were finally understood by Costabel and Dauge~\cite{Costabel-Dauge02}: in  domains with reentrant corners or edges this method may produce discrete solutions
that converge to a function that is not the solution of Maxwell's equations.  In more detail, the above problem requires us to consider convergence in the norm
(recall that $\epsilon_r$ is constant)
\[
\Vert \bfu\Vert_{\bfV_h^{(1)}}^2:=
\Vert \curl\bfu\Vert_{\Omega}^2+
\Vert \divv\bfu\Vert_{\Omega}^2+
\Vert \bfu\Vert_{\Omega}^2.
\]
Necessarily the solution $\bfE_h$ of (\ref{maxweak1}) will have
$\Vert \bfE_h\Vert_{\bfV_h^{(1)}}^2<\infty$ and hence lie in the space
 $\bfH:=\{\bfu\in \bfH_0({\rm curl};\Omega)\;|\; \Vert\bfu\Vert_{\bfV_h^{(1)}}<\infty\} $. In a convex domain, $\bfH=\bfV$ and convergence to the true solution $\bfE$ will occur.  But if $\Omega$ has reentrant corners or edges then $\bfH$ is a closed subspace of $\bfV$ with infinite codimension~\cite{Costabel-Dauge02}.  Thus there are solutions of \cref{Maxwell} that cannot be approximated by finite elements in $\bfV_h^{(1)}\subset \bfH$.  Worse, in this bad case, the finite element solutions will converge as $h$ decreases, but to a function that 
 does not solve 
Maxwell's equations.  Understanding this subtle problem was a major contribution by mathematics
to computational electromagnetism.

Of course it is possible to ``rehabilitate'' continuous elements.  Indeed
Costabel and Dauge~\cite{Costabel-Dauge02}  achieve this by weakening the divergence penalty term replacing $\gamma(\divv\bfE_h,\divv \bfxi)$ by
 $(\tilde{\gamma}\divv\bfE_h,\divv \bfxi)$ where now $\tilde{\gamma}\geq 0$ is a position dependent weight function that tends to zero approaching a reentrant corner or edge in the right way~\cite{Costabel-Dauge02}. 
 
 An alternative to using weighted spaces is instead to use negative norms of the divergence to stabilize the system~\cite{Bonito-Guermond11}.  Another alternative approach is to supplement the $H^1$ approximation space using singular functions~\cite{Ass-Cia-Son-1998}.
 
 As we have seen, across surfaces where $\epsilon_r$ is discontinuous, the normal component of the electric field will generally also jump.  So to handle discontinuous $\epsilon_r$ with $H^1$ elements it is necessary to break the continuity of the finite element functions across
 material interfaces.  Then suitable transmission
 conditions can be imposed by Nitsche's method~\cite{Assous12} which is essentially the same as the IPDG method discussed in Section~\ref{IPDG}.

\subsection{Edge Elements}\label{edgel}
In the previous section we described the use of standard continuous piecewise polynomial finite 
elements, and outlined some of the issues associated with them.  While it has proved possible to obtain an appropriate variational
formulation, and make technical modifications to
enable the use of these elements even when $\epsilon_r$ is discontinuous, some of the simplicity is lost.  In this subsection we describe the standard conforming finite element approach to approximating the time harmonic Maxwell system but which drops the goal of working in $\bfV$ and instead seeks conforming elements in the larger space
$\bfH_0({\rm curl};\Omega)$.

Finite elements in $\bfH({\rm curl};\Omega)$ date back to the work of N\'ed\'elec~\cite{ned80,ned86}.  In \cite{ned80}
he proposed a family of ``first kind'' elements 
that are obtained by augmenting function in
$\mathbb{P}_{k-1}^3$ by some extra polynomials
in $\mathbb{P}_k^3$ to give a conforming family of elements in $\bfH({\rm curl};\Omega)$.  The lowest
order finite element space (when $k=1$) is given by
\[
\bfV^{N^1_{1}}_h=\{\bfu\in \bfH_0({\rm curl};\Omega)\;|\; \bfu|_K=\bfa_K+\bfb_K\times\bfx
\mbox{ for }\bfa_K,\bfb_K\in \mathbb{C}^3\mbox{ and all }
K\in\cT_h\}
\]
Functions of this type also appear in the work of Whitney~\cite{whitney57}, and so these elements
are sometimes called Whitney elements, N\'ed\'elec
elements or edge elements.  This  family is likely the most commonly used
choice of elements, although usually
higher order spaces are used in practice~\cite{Monk03}. An important work on the mathematical foundations of
edge elements and the use of these elements for Maxwell's equations is due to Hiptmair~\cite{hip2000}, and for a different approach see~\cite{Monk03}.

Because of their relevance to discontinuous Galerkin methods, in this paper we shall consider the second family of edge elements of degree $k$ proposed by N\'ed\'elec in 1986~\cite{ned86}:
\[
\bfV^{N^2_{k}}_h=\{\bfu\in \bfH_0({\rm curl};\Omega)\;|\; \bfu|_K\in \mathbb{P}_k^3\mbox{ and all }
K\in\cT_h\}.
\]
This family differs from the first family by the 
addition of the gradients of certain polynomials.
The degrees of freedom (DOF) for this element are associated with edges, faces and tetrahedra in the mesh.  In particular, specification of the
following quantities uniquely determines a polynomial in $\mathbb{P}_k^3$ and globally guarantees a conforming space~\cite{Monk03}:
\begin{subequations}\label{edgedof}
\begin{eqnarray}
&&\Big\{\int_e(\bfu\cdot\bftau)q\,ds\mbox{ for all } q\in \mathbb{P}_k(e)
\mbox{ and all edges }e\mbox{ of }K\Big\},\label{cct2dof1}\\
&&\Big\{\int_F\bfu_T\cdot\bfq\,d A\mbox{ for all } \bfq\in D_{k-1}(F)
\mbox{ and all faces }F\mbox{ of }K\Big\},\label{cct2dof2}\\
&&\Big\{\int_K\bfu\cdot\bfq\,dV\mbox{ for all } \bfq
\in D_{k-2}(K)\Big\},
\label{cct2dof3}
\end{eqnarray}
\end{subequations}
where $\bftau$ is a unit tangent vector to the edge $e$. The spaces $D_{k-1}(F)$ and $D_{k-2}(K)$ are defined using the space of homogeneous polynomials of degree $k$ (denoted $\tilde{\mathbb{P}}_k$) as follows: $D_{k-1}(F)=(\mathbb{P}_{k-2}(F))^2\oplus \tilde{\mathbb{P}}_{k-1}(F)\,\bfx$ of 
vector functions tangential to $F$, and
$D_{k-1}(K)=(\mathbb{P}_{k-2}(K))^3\oplus \tilde{\mathbb{P}}_{k-1}(K)\,\bfx$ (see \cite{Monk03}).  

We remark that the DOFs in (\ref{cct2dof1}) are associated with the edges of the mesh and this accounts for the name ``edge elements''. Although the DOFs defined in \cref{edgedof} implicitly define a basis for $\bfV_h^{N_k^2}$ they are not convenient (at least for $k>2$) and other basis functions more suited to implementation have been defined (for example those in \cite{zaglmayr}).

Taken together, for any sufficiently smooth vector function $\bfu$ on $\Omega$, the degrees of freedom (\ref{edgedof}) define an interpolation operator $\bfr_h$ such that $\bfr_h\bfu\in \bfV_h^{N_k^2}$.  Unfortunately the interpolant is not defined for all functions in $\bfH_0({\rm curl};\Omega)$ although analogues
of the Scott-Zhang interpolation operator have been derived~\cite{CW08}.

We can now state the standard edge element
discretization of Maxwell's equations: we seek
$\bfE_h\in \bfV_h^{N^2_k}$ such that
 such that
\begin{equation}
(\curl\bfE_h,\curl\bfxi)_\Omega-\kappa^2(\epsilon_r\bfE_h,\bfxi)_\Omega=(\bfF,\bfxi)_\Omega
\mbox{ for all }\bfxi\in \bfV_h^{N^2_k}.\label{FEmaxweak}
\end{equation}
As we saw in the introduction, a key point for the variational formulation of Maxwell's equations in (\ref{maxweak})  is that the lower order term on the left hand side controls the divergence of the solution.  Edge elements mirror this at the discrete level and contain a large space of gradients.  Let
\[
S_h^{k+1}=\{p_h\in H^1_0(\Omega)\;|\;p_h|_K\in \mathbb{P}_{k+1}\}
\]
denote the space of continuous piecewise $k+1$ degree polynomials.  Since $\curl(\nabla p)=0$, it is clear that $\nabla S_h^{k+1}\subset \bfV_h^{N_k^2}$.  Choosing $\bfxi=\nabla p_h$ for any $p_h\in S_h^{k+1}$ in \cref{FEmaxweak} shows that
\[
-\kappa^2(\epsilon_r\bfE_h,\nabla p_h)_\Omega=0.
\]
Comparing to \cref{divfree}, we say that
$\epsilon_r\bfE_h$ is \emph{discrete divergence free}.  We can also write a discrete Helmholtz decomposition
\begin{equation}
\bfV_h^{N^2_k}=\tilde{\bfV}_h^{N^2_k}\oplus \nabla S_h^{k+1}
\label{discreteHelmholtz}
\end{equation}
where 
\[
\tilde{\bfV}_h^{N^2_k}=\{\bfv_h\in \bfV_h^{N_k^2}\;|\; (\epsilon_r\bfv_h,\nabla p_h)_\Omega=0\mbox{ for all }p_h\in S_h^{k+1}\}.
\]
Unfortunately functions in this space are not divergence free. 

A key property of this space is \emph{discrete compactness}. In particular, if $\{h_n\}$, $n=1,\cdots$ with $h_n\to 0$ as $n\to\infty$ and if $\bfu_n\in \tilde{\bfV}_{h_n}^{N^2_k}$ is bounded in $\bfH({\rm curl};\Omega)$, then $\{\bfu_n\}_{n=1}^\infty$ contains a subsequence converging to $\bfu\in \bfL^2(\Omega)$
and $\divv(\epsilon_r\bfu)=0$.
This property was introduced by Kikuchi~\cite{kik89} to study eigenvalue problems, and is central also to the study of the source problem.  Boffi has shown that it is equivalent to the Fortin condition from the theory of mixed methods~\cite{Boffi2000,Boffi2010}.  

Even though functions in $\tilde{\bfV}_h^{N^2_k}$ are not divergence free, if $\epsilon_r=1$,  given $\bfw_h\in \tilde{\bfV}_h^{N^2_k}$,
and if we define $\bfH\bfw_h\in \bfV$
by requiring that $\curl\bfH\bfw_h=\curl\bfw_h$ in $\Omega$ then
\begin{equation}
\Vert \bfw_h-\bfH\bfw_h\Vert_{\Omega}\leq C h^\sigma\Vert \curl\bfw_h\Vert_\Omega
\label{Herr}
\end{equation}
for some $\sigma>0$ depending on the domain $\Omega$.  Thus the discrete
divergence free function $\bfw_h$ is close to a divergence free function.
This verifies that the crucial GAP property of~\cite{Buffa} is satisfied.  Then applying the general theory developed by Buffa in~\cite{Buffa} gives the following theorem:
\begin{theorem}[Theorem 3.7 \cite{Buffa}] Suppose that $\epsilon_r$
is piecewise constant with positive real part, and $\kappa$ is not a Maxwell eigenvalue. Then for all $h$ sufficiently small there exists a unique solution $\bfE_h$ to \cref{FEmaxweak}.  In addition
\[
\Vert \bfE-\bfE_h\Vert_{{\rm curl},\Omega}\leq
C \Vert \bfE-\bfv\Vert_{{\rm curl},\Omega}
\]
for any $\bfv\in \bfV_h^{N^2_k}$.
\end{theorem}
Buffa's theory has been applied to more general scattering type problems, for example, in \cite{Gatica_meddahi}.

We have already noted that $\nabla S_h^{k+1}\subset \bfV_h^{N^2_k} $, but there is a deeper connection.  Using the following DOFs $\mathbb{P}_{k+1}$ on an element $K$:
\begin{enumerate}
\item  vertex degrees: 
$ p(\bfa_i), \quad 1 \le i \le 4$,
for the four vertices $\bfa_i$ of $K$,
\item  edge degrees: $\left\{ \frac{1}{{\rm length} (e)} \int_e p\,q \,d s\quad
\mbox{ for all } q \in P_{k - 1} (e), \; \mbox{ for all edges}\; e \mbox{ of }K \right\},
$
\item  face degrees: $ \left\{ \frac{1}{{\rm area} (f)} \int_f p\,q\, d 
A\quad \mbox{ for all }
q \in P_{k - 2} (f), \mbox{ for all faces} \; f \mbox{ of }K\right\}$,
\item  volume degrees:$
\left\{ \frac{1}{{\rm volume}(K)} \int_K p\,q \,dV\quad \mbox{ for all } q
\in P_{k - 3}  \right\},
$
\end{enumerate}
we can define an interpolation operator $\pi_h$. Then for sufficiently smooth
scalar functions $p$
\[
\nabla(\pi_h p)=\bfr_h(\nabla p).
\]
The recognition by Bossavit~\cite{Bossavit} (see also \cite{hip2000}) that this, and a similar equality for $\bfH({\rm curl};\Omega)$ and $\bfH({\rm div};\Omega)$ interpolation, constitute a discrete de Rham complex motivated in part the Finite Element Exterior Calculus of
Arnold, Falk and \rev{Winther}~\cite{afw00,afw10}. This has resulted\rev{, for example,} in a much better understanding of edge elements, new elements on \rev{pyramids~\cite{nigam-phillips}}, and a general theory \rev{for} Laplace type equations \rev{(i.e. involving the Hodge Laplacian)}.

\section{Interior Penalty DG Methods}\label{IPDG}
One of the first discontinuous Galerkin methods
to be developed for Maxwell's equations followed the 
Interior Penalty Discontinuous Galerkin (IPDG)
approach~\cite{Arnold82}. An early paper explicitly including control of the divergence~\cite{Perugia02stabilizedinterior} was soon superseded by
improvements that do not require explicit divergence control. The version we shall describe here is from~\cite{Houston2005}.

To outline the derivation of the method, we apply
\cref{intp} to an element $K$ in the mesh
and using \cref{Maxwell-curl} we obtain 
for any smooth vector function $\bfxi$ on $K$
\begin{equation}
(\curl \bfE,\curl\bfxi)_K-k^2(\epsilon_r\bfE,\bfxi)_K-(\bfF,\bfxi)_K+
\int_{\partial K}\bfnu_K\times\curl\bfE\cdot\overline{\bfxi}_T\,dA=0,
\label{DGbase}
\end{equation}
where $\bfnu_K$ is the unit outward normal to $K$.
To proceed, we will add the above equality over all elements, and to do this we need more notation.  If $K^+$ and $K^-$ meet at a face $F$
we define the standard jump and average value of a piecewise smooth vector function $\bfv$ taking values $\bfv^+$ on $K^+$ and $\bfv^-$ on $K^-$
by
\[
\jmp{\bfv}:=\bfnu_{K^+}\times \bfv^++\bfnu_{K^-}\times \bfv^-,
\qquad
\avg{\bfv}=(\bfv^++\bfv^-)/2
\]
For a face $F\subset\Gamma$, $\jmp{\bfv}=\bfnu\times\bfv$ and $\avg{\bfv}=
\bfv$.

Adding \cref{DGbase} over all elements, and using the following ``DG magic formula''~\rev{(see e.g. \cite{MelenkEsterhazy12,Perugia02stabilizedinterior})}
\[\sum_K \int_{\partial K}\bfnu_K\times\curl\bfE\cdot\overline{\bfxi}_T\,dA
= 
\int_{\cF_I}\jmp{\curl\bfE}\cdot\avg{\overline{\bfxi}}-\avg{\curl\bfE}\cdot\jmp{\overline{\bfxi}}\,dA+\int_{\cF_B}\jmp{\curl\bfE}\cdot\avg{\overline{\bfxi}}\,dA
\]
where $\cF_I$ is the \rev{union} of all interior faces of
the mesh and $\cF_B$ is the \rev{union} of all boundary 
faces, using the PEC boundary condition and
the continuity of the tangential component of $\curl\bfE$ across internal faces gives
\begin{equation}
(\curl \bfE,\curl\bfxi)_{\cT_h}-\kappa^2(\epsilon_r\bfE,\bfxi)_{\cT_h}-(\bfF,\bfxi)_{\cT_h}-
\int_{\cF_I}\avg{\curl\bfE}\cdot\jmp{\overline{\bfxi}}\,dA=0.\label{IPDG1}
\end{equation}
This identity is the basic result for writing
down the upcoming IPDG method after discretization.

To discretize the problem we introduce the 
space of discontinuous vector polynomials of degree $k$
\[
\bfV^{IP}_h=\{\bfv\in \bfL^2(\Omega)\;|\;
\bfv|_K\in\mathbb{P}_k^3\mbox{ for all } K\in\cT_h\}.
\]
Then using (\ref{IPDG1}) we are led to the IPDG method
proposed in \rev{\cite{Houston2005}}:
find $\bfE_h\in \bfV_h^{IP}$ such that
\begin{equation}
a^{IP}(\bfE_h,\bfxi)=\int_{K}\bfF\cdot\overline{\bfxi}\,dV\mbox{ for all }\bfxi\in \bfV_h^{IP}
\label{IPProb}
\end{equation}
where
\begin{eqnarray*}
a^{IP}(\bfE_h,\bfxi)&=&
(\curl \bfE_h,\curl\bfxi)_{\cT_h}-k^2(\epsilon_r\bfE_h,\bfxi)_{\cT_h}\\&&-
\int_{\cF_I}\left(\avg{\curl\bfE_h}\cdot\jmp{\overline{\bfxi}}+\avg{\overline{\curl\bfxi}}\cdot\jmp{\overline{\bfE_h}}\right)\,dA+
\int_{\cF_I}\frac{\alpha}{\rm h}\jmp{\bfE_h}\cdot\jmp{\overline{\bfxi}}\,dA.
\end{eqnarray*}
Here, following the usual IPDG philosophy~\cite{Arnold82}, we have added
a symmetrizing (but consistent) term and a term that penalizes the jumps across element faces.  In addition, $\alpha>0$ is a penalty parameter that must be chosen sufficiently large and 
\[
{\rm h}(\bfx)=h_F\mbox{ for }\bfx\in F
\]
where $h_F$ denotes the diameter of the face $F$. 

From the definition of $a^{IP}(.,.)$, the appropriate norm to measure the error in the 
solution is the DG norm
\[
\Vert  \bfv \Vert^2_{DG}=\Vert {\rm h}^{-1/2}\jmp{\bfv}\Vert_{\cF_I}^2+\Vert \bfv \Vert_\Omega^2+\sum_{K\in\cT_h}\Vert \curl\bfv\Vert_K^2{\color{red}.}
\]
Then under the condition that $\epsilon_r$ is constant in $\Omega$, the following theorem can be shown:
\begin{theorem}[\rev{Th. 3.2 of \cite{Houston2005}}] Assume that the true solution satisfies $\bfE\in \bfH^s(\Omega)$ and $\curl\bfE\in H^s(\Omega)$ for some $s>1/2$ and assume that the \rev{penalty} parameter is chosen sufficiently large.  Then for all $h$ small enough, there is a unique solution
to \cref{IPProb}.  In addition the following
error estimate holds:
\[
\Vert\bfE-\bfE_h\Vert_{DG}\leq Ch^{\min(s,k)}\left(
\Vert\bfE\Vert_{\bfH^s(\Omega)}+\Vert\curl\bfE\Vert_{\bfH^s(\Omega)}\right)
\]
\end{theorem}
While we won't give full details (see \cite{Houston2005}), the proof of this result proceeds as follows.  First, by using the definition of the sesquilinear form, the following estimate holds:
\[
\Vert\bfE-\bfE_h\Vert_{DG}\leq C\left(
\inf_{\bfv\in \bfV_h^{IP}}\Vert\bfE-\bfv\Vert_{DG}+
\mathcal{R}_h(\bfE)+\mathcal{E}_h(\bfE-\bfE_h)\right)
\]
where, after some calculation, and defining $\boldsymbol{\Pi}_h$ to be the $\bfL^2$ projection onto $\bfV_h^{IP}$, we have (see Lemma \rev{4.9 of \cite{Houston2005}})
\[
\mathcal{R}_h(\bfE)=\sup_{\bfv\in \bfV_h^{IP} }
\frac{\int_{\cF}\jmp{\bfv}\cdot\avg{\curl\bfE-\boldsymbol{\Pi}_h\curl\bfE}\,dA}{\Vert \bfv\Vert_{DG}}.
\]
This can then be estimated in the usual way using error estimates for $\boldsymbol{\Pi}_h$.

The term $\mathcal{E}_h(\bfE-\bfE_h)$ is given by
\[
\mathcal{E}_h(\bfE-\bfE_h)=\sup_{\bfv\in \bfV_h^{IP}}
\frac{|(\bfE-\bfE_h,\bfv)_\Omega|}{\Vert \bfv\Vert_{DG}}
\]
and arises from the lower order term in Maxwell's
equations (as is usual for time harmonic wave equations). In \rev{\cite{Houston2005}}
this is analyzed using an adjoint problem (see \cite{Monk03} for a similar analysis in the case of edge elements).  In particular, for a given $\bfv\in\rev{\bfV}_h^{IP}$, let $\bfv^c \in \bfV_h^{N_k^2}$ be a conforming edge element approximation of $\bfv$ such that (see Prop. \rev{4.6 of \cite{Houston2005}})
\begin{equation}
\Vert \bfv -\bfv^c\Vert_{\Omega}\leq Ch\Vert \bfv\Vert_{DG}.\label{DGerr}
\end{equation}
Then let $\bfv^c=\bfv^c_0+\nabla p_h$ according to the discrete Helmholtz decomposition \cref{discreteHelmholtz} (i.e. $\bfv_0^c$ is discrete divergence free).  Using the operator $\bfH$ from 
\cref{FEMax} we have (cancelling the gradient term because $\bfE-\bfE_h$ is discrete divergence free)
\[
(\bfE-\bfE_h,\bfv)=(\bfE-\bfE_h,\bfv-\bfv^c)+(\bfE-\bfE_h,\bfv_0^c-\bfH\bfv_0^c)+(\bfE-\bfE_h,\bfH\bfv_0^c).
\]
The first term on the left hand side is approximated
by \cref{DGerr}, the second from the \cref{Herr} and the third is estimated using duality via the solution $\bfz\in H_0({\rm curl};\Omega)$ of the adjoint Maxwell problem
\[
\curl\curl\bfz-\kappa^2\overline{\epsilon_r}\bfz=\bfH\bfv_0^c\mbox{ in }\Omega.
\]
The key here is that the right hand side is exactly divergence free, and so $\bfz$ is smoother than $\bfE-\bfE_h$.  Using the adjoint problem, and approximation properties of the N\'ed\'elec interpolant $\bfr_h$ of $\bfz$, it is then possible to show that
\[
|(\bfE-\bfE_h,\bfH\bfv_0^c)|\leq C h^\sigma
\Vert\bfE-\bfE_h\Vert_{DG}\Vert\bfv\Vert_{DG}
\]
for some $\sigma>0$ and the usual kickback argument completes the proof.

The direct use of the  adjoint problem here \rev{limits the proof to constant $\epsilon_r$. For the case of piecewise smooth  $\epsilon_r$, a more sophisticated analysis in \cite[Section 6]{Buffa+Perugia}  proves an inf-sup
condition in this case, and results in a quasi-optimal error estimate (see in particular  \cite[Remark 7.11]{Buffa+Perugia}).}  Our discussion was intended  illustrates the interplay of 
conforming and DG approximations and the use of duality.

\section{HDG Methods}\label{HDG}
In the Introduction we have argued that DG methods offer several advantages compared to conforming methods. However, the total number of global degrees of freedom (DOF)
of a DG method is much more than that of an $\bm H({\rm curl};\Omega)$-conforming edge element method on the same mesh and of the same order.

\rev{Hybridizable discontinuous Galerkin (HDG) methods \cite{Cockburn_Gopalakrishnan_Lazarov_Unify_SINUM_2009} were recently
introduced with the aim of reducing the dimension of the discrete global linear system that needs to be solved and deliver superconvergent recovery of variables of interest.} By design, degrees of freedom associated with the volume of each element can be eliminated from the global discrete problem in a process akin to static condensation. The resulting global system of an HDG method only involves DOFs on the skeleton of the mesh (i.e. all faces of the mesh) and so the global system after condensation is much smaller.

\rev{We now follow closely \cite{Cockburn_maxwell_HDG_JCP_2011} in order} to write down a standard HDG method for Maxwell's equations. We first introduce a vector variable $\bft=\nabla\times\bm E$ (this is just a scaled magnetic field variable which is a quantity of interest in most simulations). We can then rewrite \eqref{Maxwell} into the following first order system: find  $(\bft,\bm E)$ such that
\begin{subequations}\label{mixed}
	\begin{align}
		\bft-\nabla\times\bm E&=\bm 0&\text{in }\Omega,\\
		\nabla\times\bft-\kappa^2\epsilon_r\bm E&=\bm F&\text{in }\Omega,\\
		 \bfnu\times\bm E&=\bm 0&\text{on }\partial\Omega.
	\end{align}
\end{subequations}

Given a choice of three finite dimensional polynomial spaces $\bm V(K)\subset \bm H^1(K)$, $\bm W(K)\subset \bm H(\text{curl}; K)$ and $ \bm M(F)\subset \bm L^2(F)$,  where $K$ is an arbitrary element in the mesh and $F$ is an arbitrary \rev{face (or edge in 2D)}, we define the global spaces by
\begin{align*}
	\bm V_h&:=\{ \bm v\in \bm L^2(\mathcal{T}_h): \bm v|_K\in \bm V(K), K\in\mathcal{T}_h\},\\
	\bm W_h&:=\{ \bm w\in  \bm L^2(\mathcal{T}_h): \bm w|_K\in \bm W(K), K\in\mathcal{T}_h\},\\
	\bm M_h&:=\{\bm \mu\in\bm L^2(\mathcal E_h): \bm \mu|_F\in\bm M(F), F\in\mathcal{F}_h\}.
\end{align*}
Note that if \rev{$\bfW(K)=[\mathbb{P}_k(K)]^d$ where $d=2,3$ is the spatial dimension,} then $\bfW_h=\bfV_h^{IP}$ introduced earlier.

We can now derive the HDG method for (\ref{mixed}) by multiplying each equation by the complex conjugate of an appropriate discrete test
function, integrating element by element and using the integration by parts identity \cref{intp} element by element in the usual way. Summing the results over all elements, the  HDG methods seeks an approximation to $(\bft,\bm E, \bm E|_{\mathcal{F}_h})$, by $(\bft_h,\bm E_h, \widehat{\bm E}_h)\in \bm V_h\times\bm W_h\times  \bm M_h$, such that
\begin{subequations}\label{Maxwell_equation_HDG_form_ori}
	\begin{align}
		(\bft_h, \bm v_h)_{\mathcal{T}_h}-(\bm E_h,\nabla\times \bm v_h)_{\mathcal{T}_h}-\langle \bm n\times\widehat{\bm E}_h, \bm v_h \rangle_{\partial\mathcal{T}_h}&=0,\label{Maxwell_equation_HDG_form_ori_1}\\
		(\bft_h,\nabla\times\bm w_h)_{\mathcal{T}_h}
		+\langle\bm n\times\widehat{\bft}_h,{\bm w}_h \rangle_{\partial\mathcal{T}_h}-(\kappa^2\epsilon_r\bm E_h,\bm w_h)_{\mathcal{T}_h}&=(\bm F,\bm w_h)_{\mathcal{T}_h},\label{Maxwell_equation_HDG_form_ori_2}\\
		\langle\bfnu\times\widehat{ \bft}_h,\widehat{\bm w}_h \rangle_{\partial \mathcal T_h/\partial\Omega}&=0,\label{Maxwell_equation_HDG_form_ori_3}\\
		\langle\bfnu\times\widehat{\bm E}_h,\bm n\times\widehat{\bm w}_h \rangle_{\partial\Omega}&=0\label{Maxwell_equation_HDG_form_ori_4}
	\end{align}
	for all $(\bfr_h,\bm v_h, \widehat{\bm v}_h)\in \bm V_h\times\bm W_h\times\bm M_h$. 
\end{subequations}	

Different choices of the space $\bm V(K)\times \bm W(K)\times \bm M(F)$ and numerical flux  $\bfnu\times\widehat {\bft}_h$ give different HDG methods.

Two HDG methods were presented in \cite{Cockburn_maxwell_HDG_JCP_2011} to approximate  Maxwell's equations. Both these HDG methods use 
\begin{align}\label{space_specify}
   \bm V(K)\times \bm W(K)\times \bm M(\rev{F}) = [\mathbb{P}_k(K)]^d\times [ \mathbb{P}_k(K)]^d\times [ \mathbb{P}_k(F)]^{d-1},
\end{align}
where $d=2,3$ is the space dimension.  The numerical flux is written including a stabilization term as
\begin{align}
	\bfnu\times\widehat{\bft}_h=
	\bfnu\times \bft_h+\tau(\bm E_h-\widehat{\bm E}_h)_T,\label{Maxwell_equation_HDG_form_ori_5}
\end{align}
where $\tau$ is a positive constant.

The first HDG method of \cite{ Cockburn_maxwell_HDG_JCP_2011} extends (\ref{Maxwell_equation_HDG_form_ori}) by enforcing the divergence-free condition on the electric induction $\epsilon_r\bfE$ and introduces a Lagrange multiplier to accomplish this. After hybridization (static condensation), it produces a linear
system for the DOF of the approximate traces of both the
tangential components of the vector field $\hat{\bfE}_h$ and the global Lagrange multiplier. The second
HDG method does not enforce the divergence-free condition and results in a linear
system only for the DOF of the approximate trace of the tangential components of
the vector field $\hat{\bfE}$. The well-posedness, conservativity and consistency of the two HDG methods, together with a numerical demonstration, was shown in \cite{Cockburn_maxwell_HDG_JCP_2011}. However, this paper does not include an error analysis and the numerical experiments are only for $d=2$ (two dimensional domains). 

The case of HDG on a three dimensional domain was discussed in \cite{MR3117424}.  The resulting linear system was solved using a domain decomposition technique but without a convergence analysis.

There are some papers with a complete convergence analysis. The first was proposed in \cite{Feng_Maxwell_CMAM_2016}: the authors use $k=1$, $d=3$ and $\tau=1/h$ in \eqref{space_specify} and \eqref{Maxwell_equation_HDG_form_ori_5}, respectively. A detailed $hp$ and a posteriori convergence analysis is found in  \cite{Lu_hp_HDG_Maxwell_Math_Comp_2017,Chen_Maxwell_HDG_CMAME_2018} for arbitrary polynomial degree and $d=3$.  However, these works only obtained a suboptimal convergence rate for $\bft$. 

In a very recent paper~\cite{ChenMonkZhang1}, we used the  concept of an $M$-decomposition, which was proposed by Cockburn et al in \cite{Cockburn_M_decomposition_Part1_Math_Comp_2017} for elliptic PDEs \rev{to analyze HDG schemes
for Maxwell's equations in two dimensions.  This analysis provides conditions on the HDG spaces need to obtain optimal convergence, and superconvergence of some variables.  The extension of this approach to 3D is challenging, and remains to be done.}

Another way to design an  optimally convergent HDG method is  to  define a  special numerical flux and the new spaces 
$\bm V(K)$, $\bm W(K)$ and $\bm M(\partial K)$, even though the spaces do not admit $M$-decompositions. 
In \cite{Chen_maxwell_HDG_2017}, Chen et al proposed a new HDG method for Maxwell's equations by augmenting the boundary space and choosing
\begin{align*}
   \bm V(K)\times \bm W(K)\times \bm M(\partial K) =  [\mathbb{P}_k(K)]^3\times [ \mathbb{P}_{k+1}(K)]^3\times [\mathbb{P}_k(F)\oplus\nabla \widetilde{\mathbb{P}}_{k+2}(F)]^2,
\end{align*}
where 
$\widetilde{\mathbb{P}}_{k+2}(F)$ denotes the set of homogeneous
polynomials of degree $k + 2$ on $F$. The stabilization has the form 	
\begin{align}
	\bm n\times\widehat{\bft}_h=
	\bm n\times \bft_h+\frac{1}{\bm h}(\bm P_M\bm E_h-\widehat{\bm E}_h      )_T,\label{Maxwell_equation_HDG_form_ori_6}
\end{align}
where $\bm P_M$ is an elementwise $L^2$ orthogonal projection into $\bm M(F)$. They obtain an optimal convergent rate for the solution by some regularity assumptions on the dual problem at the expense of a larger trace space.

\section{The Ultra Weak Variational Formulation}\label{UWVF}
Discontinuous Galerkin methods are not limited to piecewise polynomial elements.  In this section we present a DG method that uses
solutions of Maxwell's equations element by element.  Methods that use solutions of the underlying equation in the approximation scheme are termed Trefftz methods~\cite{Trefftz26}, and the method we shall describe falls into the general class of Trefftz DG schemes.  We shall also consider a more general version of the boundary value problem that brings us closer to practical applications.

The particular method we will describe is the Ultra Weak Variational Formulation (UWVF) of
Maxwell's equations due to Cessenat~\cite{cessenat_phd} and developed
further in \cite{HMM06}. 
This method applies to a more general problem that includes the previously
considered PEC boundary condition as a special case, but is restricted to equations in which $\bfF=0$. In addition it is necessary to assume that $\epsilon_r$ is piecewise constant on the mesh.  Therefore we now assume that $\bfE$ satisfies
\begin{subequations}\label{Maxwell-UWVF}
\begin{align}
\curl\curl\bfE-\kappa^2\epsilon_r
\bfE&=0\mbox{ in }\Omega,\label{Maxwell-curl-UWVF}\\
\bfnu\times\curl\bfE-i\kappa\lambda \bfE_T&=Q(\bfnu\times\curl\bfE+i\kappa\lambda \bfE_T)+\bfg\mbox{ on }\Gamma,\label{Maxwell-UWVF-BC}
\end{align}
\end{subequations}
where $\lambda>0$ is a real parameter, $Q$ is another real parameter with $|Q|<1$ and
$\bfg\in \bfL^2_T(\Gamma):=\{\bfg\in \bfL^2(\Gamma)\;|\; \bfnu\cdot\bfg=0\}
$ is data.  Note that $Q=1$ gives the PEC boundary condition considered before.

To derive the UWVF
we apply \cref{intp} again to the curl term in
\cref{DGbase} (with $\bfF=0$) to obtain
\begin{equation}
(\bfE,\curl\curl\bfxi-\kappa^2\overline{\epsilon_r}\bfxi)_K+
\int_{\partial K}\bfnu_K\times\curl\bfE\cdot\overline{\bfxi}_T+\bfnu_K\times\bfE\cdot\overline{\curl\bfxi}_T\,dA=0.
\label{UVWFbase}
\end{equation}
Clearly, if $\bfxi$ solves the adjoint
Maxwell system
\begin{equation}
\curl\curl\bfxi-\kappa^2\overline{\epsilon_r}\bfxi=0\mbox{ in }K,\label{maxadj}
\end{equation}
then we obtain the identity
\[
\int_{\partial K}\bfnu_K\times\curl\bfE\cdot\overline{\bfxi}_T+\bfnu_K\times\bfE\cdot\overline{\curl\bfxi}_T\,dA=0.
\]
Using this identity we can prove the following ``isometry'' result (see \cite{cessenat_phd}):
\begin{eqnarray}
&&\int_{\partial K}\frac{1}{\lambda}\left(
\bfnu_K\times\curl\bfE+i\kappa\lambda \bfE_T\right)\cdot\left(\overline{\bfnu_K\times\curl\bfxi+i\kappa\lambda \bfxi_T}\right)\,dA\nonumber\\&
=&
\int_{\partial K}\frac{1}{\lambda}\left(
\bfnu_K\times\curl\bfE-i\kappa\lambda \bfE_T\right)\cdot\left(\overline{\bfnu_K\times\curl\bfxi-i\kappa\lambda \bfxi_T}\right)\,dA.\label{iso}
\end{eqnarray}
Now consider an arbitrary element $K$ in the mesh. The four faces are either faces of another element in the mesh (interior faces) or boundary faces.  Using the above equality and taking into account the difference in sign of the normal vector on adjacent elements we obtain
\begin{eqnarray*}
&&\int_{\partial K}\frac{1}{\lambda}\left(
\bfnu_K\times\curl\bfE+i\kappa\lambda \bfE_T\right)\cdot\left(\overline{\bfnu_K\times\curl\bfxi+i\kappa\lambda \bfxi_T}\right)\,dA\\&
=&
\int_{\partial K\cap \partial K'}\frac{1}{\lambda}\left(-
\bfnu_{K'}\times\curl\bfE-i\kappa\lambda \bfE_T\right)\cdot\left(\overline{\bfnu_{K}\times\curl\bfxi-i\kappa\lambda \bfxi_T}\right)\,dA\\
&&+\int_{\partial K\cap \Gamma}\frac{1}{\lambda}\left[Q\left(
\bfnu_K\times\curl\bfE+i\kappa\lambda \bfE_T\right)+\bfg\right]\cdot\left(\overline{\bfnu_K\times\curl\bfxi-i\kappa\lambda \bfxi_T}\right)\,dA.
\end{eqnarray*}
Assuming that $\bfE$ is sufficiently regular such that $\bfnu\times\curl\bfE+i\kappa\lambda \bfE_T\in \bfL_T^2(\partial K)$, and defining the space
\[
\bfW=\Pi_{K\in\cT_h}\bfL_T^2(\partial K),
\]
we define the unknown boundary \rev{impedance} flux $\bfchi_K=\bfnu\times\curl\bfE|_{K}+i\kappa\lambda (\bfE|_K)_T$ so that $\bfxi:=\Pi_{K\in \cT_h}\bfxi_K\in \bfW$
satisfies
\begin{eqnarray}
&&\int_{\partial K}\frac{1}{\lambda}\bfchi_K\cdot\left(\overline{\bfnu_K\times\curl\bfxi+i\kappa\lambda \bfxi_T}\right)\,dA\nonumber\\&
=&
-\int_{\partial K\cap \partial K'}\frac{1}{\lambda}\bfchi_{K'}\cdot\left(\overline{\bfnu_{K}\times\curl\bfxi-i\kappa\lambda \bfxi_T}\right)\,dA\nonumber\\
&&+\int_{\partial K\cap \Gamma}\frac{1}{\lambda}\left[Q\bfchi_K+\bfg\right]\cdot\left(\overline{\bfnu_K\times\curl\bfxi-i\kappa\lambda \bfxi_T}\right)\,dA\label{eqUWVFeq}
\end{eqnarray}
for  all test
functions $\bfxi$ such that
$\bfnu_{K}\times\curl\bfxi+i\kappa\lambda \bfxi_T\in \bfL^2_T(\partial K)$ and $\bfxi$ satisfies \cref{maxadj} on $K$, for  all $K\in\cT_h$.

The key idea from \cite{cessenat_phd} is to
use a space of plane wave solutions to the adjoint problem to discretize  \eqref{Maxwell-UWVF}.  For each element $K\in\cT_h$
we take $p_K$ independent real direction vectors $\bfd_j^K$, $|\bfd_j^K|=1$ and for each
$j$ we choose two mutually orthogonal real unit polarization
vectors $\bfp_{j,\ell}^K$ with $\bfp_{j,\ell}^K\cdot\bfd_j^K=0$, $\ell=1,2$.  Then the field
$\bfp_{j,\ell}\exp(i\kappa\sqrt{\overline{\varepsilon_r}}\bfd_j^K\cdot\bfx)$ satisfies the adjoint problem in $K$. Then we define the local plane wave space on an element $K$ by
\[
\bfPW_h^K=\mbox{span}\left\{
\bfp^K_{j,\ell}\exp(i\kappa\sqrt{\overline{\epsilon_r}}\bfd_j^K\cdot\bfx)\;|\;1\leq j\leq p_K,\; \ell=1,2\right\}.
\]
Using this space, the global test and trial space is
\[
\bfW_h=\Pi_{K\in\cT_h}
\{\bfnu\times\curl\bfxi_h+i\kappa\lambda (\bfxi_h)_T\;|\; \bfxi_h\in \bfPW^K_h\}.
\]
With this in hand, we seek $\bfchi_h\in \bfW_h$ such that
(\ref{eqUWVFeq}) holds for all $\bfxi_h\in \bfPW_h^K$ and all elements $K\in \cT_h$.

The choice of direction vectors $\bfd_j^K$ can be made in a variety of ways.  An example with good theoretical properties is given in
\cite{Hiptmair2011ErrorAO}, but we use the convenient Hammersley points on the unit sphere.

Cessenat~\cite{cessenat_phd} shows that, if $|Q|<1$, this discrete UWVF has a unique solution and derives error estimates on $
\Gamma$. A computational study relating the UWVF to upwind DG methods, and deriving the Perfectly Matched Layer (a mesh truncation technique) in this case was performed in
\cite{HMM06}.

Historically, regarding error analysis, in \cite{buf07,git07}, it was noted that the UWVF for the Helmholtz equation is a special case of the IPDG method (i.e. a special choice of parameters in IPDG) when a Trefftz plane wave basis is used and when the coefficient functions in the Helmholtz equation are real.  This leads to an error analysis for the Helmholtz equation (see for example~\cite{perugiahp} and the survey article \cite{Hiptmair2016}).  

The UWVF for Maxwell's equations is also equivalent to an IPDG scheme when $\epsilon_r$ is real (and of course, piecewise constant).
This theory depends on having the same plane wave
space for trial and test space, and the equivalence of UWVF and IPDG does not
hold for absorbing media (complex parameters).  The plane wave IPDG scheme was analyzed in \cite{Hiptmair2011ErrorAO} using the interesting stability result from \cite{moi11} and convergence was verified. Hence convergence of the UWVF is also verified. The convergence of UWVF for complex $\epsilon_r$ is not proved.

Because of the isometry result (\ref{iso}), solving the discrete UWVF system can be accomplished by simple iterative techniques~\cite{cessenat_phd,HMM06}.  Although not justified theoretically, we typically use the BiCGstab algorithm~\cite{vorst92}.
This takes many hundreds of iterations but is easy to parallelize.

In comparison to HDG, we can see that HDG is a single trace formulation  ($\hat{\bfE}_h$) whereas UWVF involves two traces $\bfchi_h$ on a face between an element $K$ and its neighbor $K'$.

The main drawback of the UWVF and the Trefftz IPDG scheme is that the discrete problem becomes rapidly ill-conditioned as $p:=\max_{K\in\cT_h}p_k$ increases.
In \cite{HMM06} we used several ad-hoc 
techniques to control conditioning by varying $P_K$ from element to element, but this in turn limits accuracy.  It is likely that a combined Trefftz and finite element basis (Trefftz on some larger 
elements and finite elements on small elements or near singularities) could be attractive~\cite{mon10}.

Modern implementations of UWVF use 
general element types (besides simplices, also prisms and hexahedra) to help with
meshing layers such as the PML.  In principle much more complicated elements
are allowed by the convergence theory~\cite{Hiptmair2011ErrorAO}.

As an example of the use of the UWVF, we consider the problem of computing exterior
scattering from a unit ball with $\epsilon_r=1$ and $\kappa=44$. The parameter $\lambda=1$ and $Q$ is set to $Q=1$ on the surface of the sphere (PEC boundary condition).  The computational domain
is the annular region inside the cube $[-1.8568,1.8568]^3$ outside of the ball.  We choose $Q=0$ 
on the outer surface of the domain (this is a simple absorbing boundary condition
that approximates scattering on an infinite domain).
In addition, a PML with parameter $\sigma=2$ is used to cut down 
reflection outside the cube $[-1.5712,1.5712]^3$ (see \cite{HMM06} for how to implement a
simple PML in the UWVF). The  mesh  consists of 127,113
tetrahedra.  The number of directions
per element is between 29 and 67 and is chosen to equilibrate the local condition number
of the inner product matrix on each element~\cite{HMM06}. The directions themselves are Hammersley points. This results in 10,743,064 degrees of freedom.
The BICG solver used 154 iterations to reduce the residual by a factor of $10^{-5}$.
This calculation was performed in 432 seconds using 20 MPI parallel processes
 on a 20 core Linux computer using an Intel(R) Xeon(R) Gold 6138 CPU at 2.00GHz.
This time includes assembling the matrices, solving via BiCGstab and computing the far field pattern at 720 points.

\begin{figure}
    \centering
    \begin{tabular}{cc}
    \resizebox{0.48\textwidth}{!}{\includegraphics{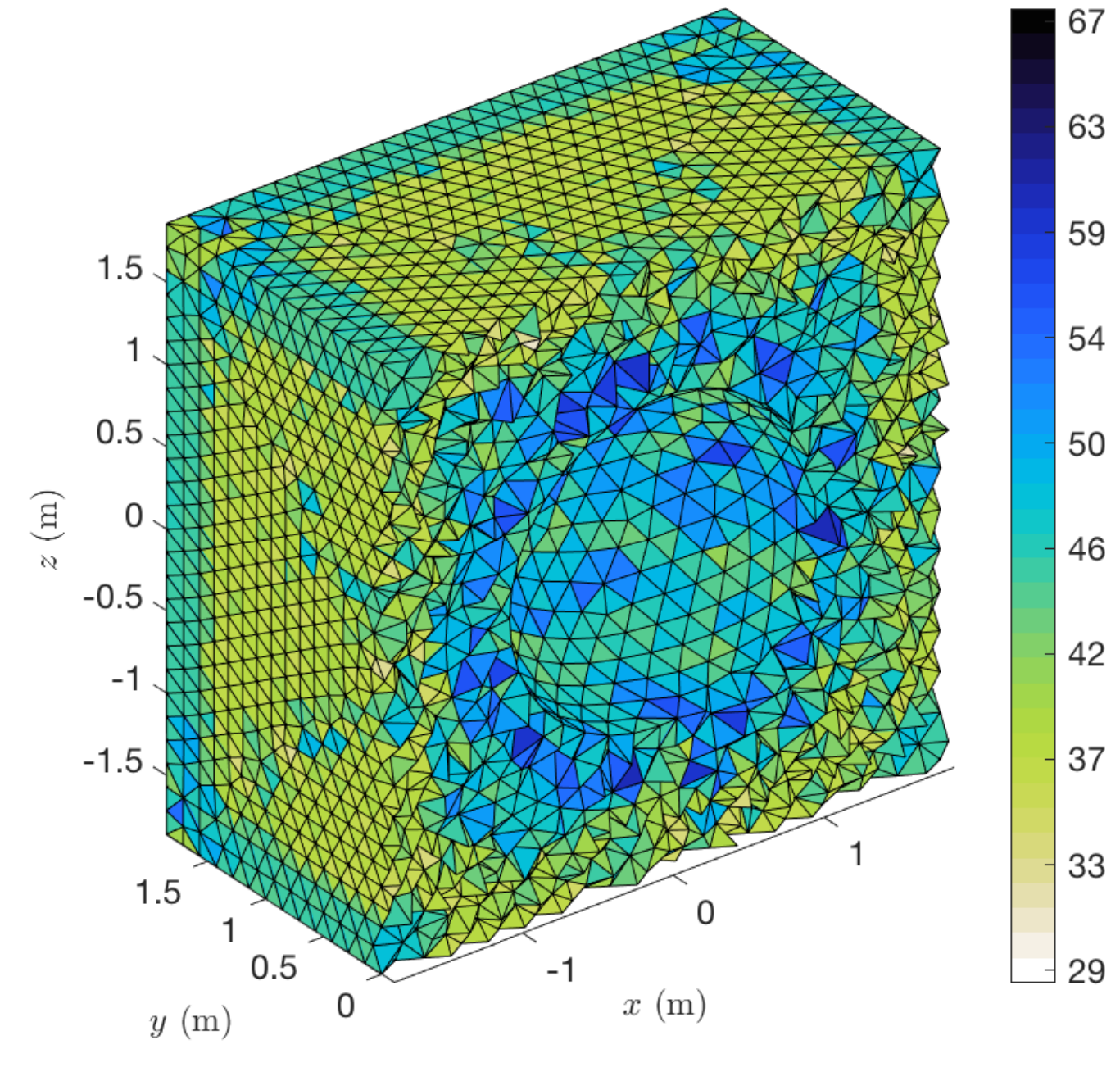}}&   \resizebox{0.48\textwidth}{!}{\includegraphics{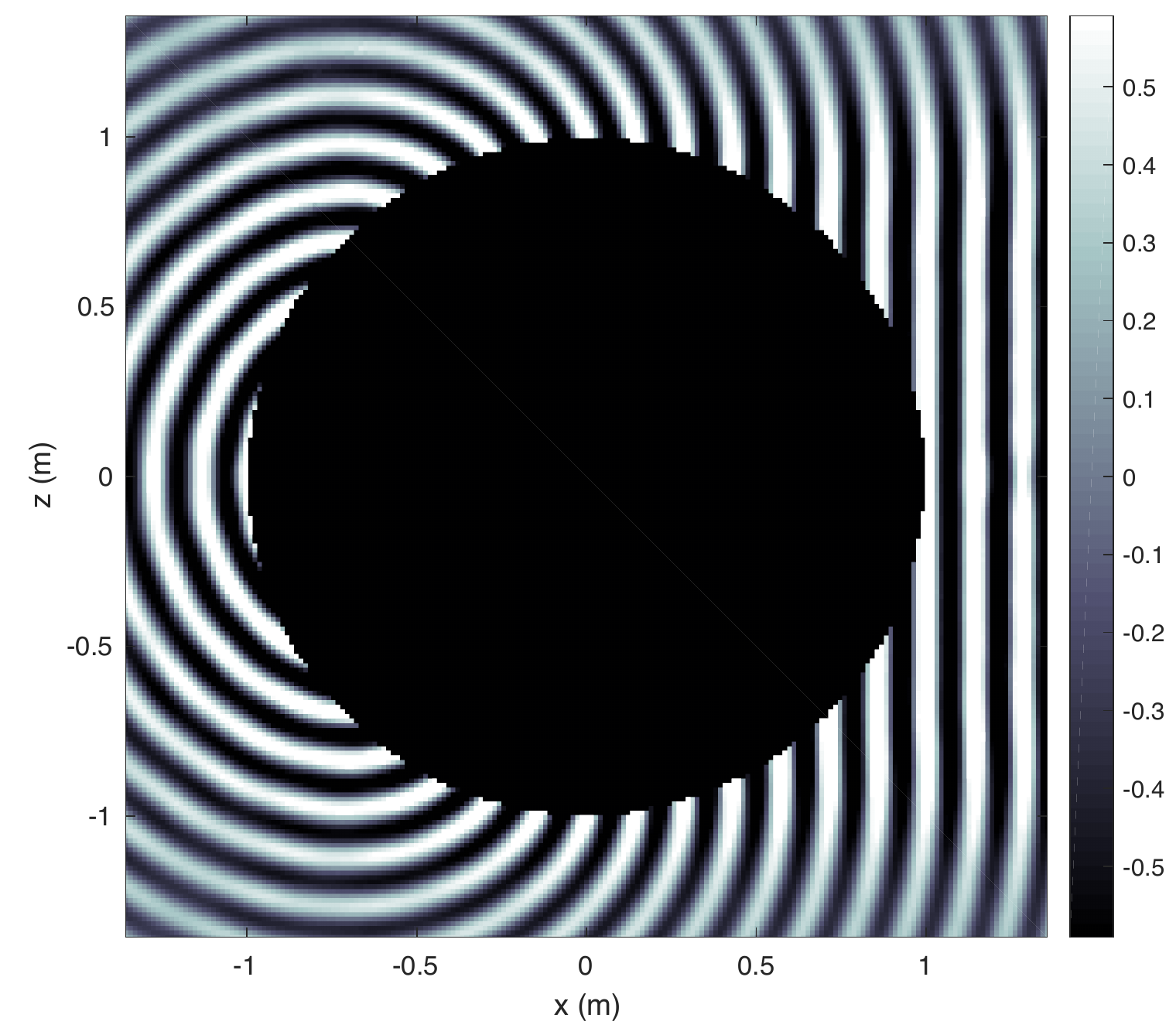}}
    \end{tabular}
    \caption{Left: A section through the mesh used in the calculation of scattering from a sphere. The shading shows the number $p_K$ of directions on each element.  Right: The real part of $E_y$ as function of $x-z$ in the plane $y=0$ (only the 
    region inside the PML is plotted).}
    \label{fig1}
\end{figure}

In Fig.~\ref{fig1}, the left panel shows a slice through the mesh. The right hand panel shows the real part of the $y$ component of the scattered electric
field $\bfE_h$ in the $x-z$ plane that was created by an incident plane wave along the $x$-axis.  This is a typical
model problem for scattering calculations
since the exact solution is known~\cite{Monk03}.
Generally, it is more efficient to use the UWVF with elements that have a geometric diameter that is at least one wavelength.  This may then introduce an unacceptable 
error due to approximating any curved boundaries of scatterers by facets.  Therefore it is essential to include curved elements in the mesh.  Although numerical integration must be used on curved faces of elements, this is simplified compared to polynomial based schemes because only surface integration needs to be performed.

The example shows some of the issues for a useful code which must be able to approximate scattering problems and also handle much more complex domains and higher frequencies than that shown here.  Of course, realistic 
calculations generally demand many more tetrahedra and degrees of freedom.

\section{Conclusions and the Future}

From the point of view of error analysis, the conforming edge finite element method is the best understood.  Furthermore, such elements provide critical conforming approximation
results for proofs of convergence for other methods such as IPDG.  However, the three DG schemes we have discussed are less well justified theoretically, especially for inhomogeneous media.  

All the analysis mentioned here does not
explicitly track the $\kappa$ dependence of
the constants appearing in the error bounds.
Obtaining $\kappa$ dependent bounds for Maxwell solvers along the lines of those for the Helmholtz equation in~\cite{Sauter_Melenk} would also be interesting
from the point of view of high wave number applications.

Looking ahead, the main problem with all these methods is to obtain a fast solver for the discrete matrix problem after discretization.  \rev{Numerical results for the UWVF (not shown) show that the number of iterations to solve the 
discrete problem are not negatively affected by changes in wave number $\kappa$ for fixed $h$ and $p$. This UWVF results suggests some hope
that other} hybridized solvers, working only on the skeleton of the mesh, may help with
the $\kappa$ dependence of the iteration number.

\section*{Acknowledgements}
The research of P.B. Monk and Y. Zhang is partially supported by the US National Science Foundation (NSF) under grant number DMS-1619904.
Research of P.B. Monk is also supported  in part by AFOSR under grant FA9550-17-1-0147.

\bibliographystyle{amsplain}
\bibliography{MZMoC75}

\end{document}